\documentclass[12pt]{amsart}

\newtheorem{theorem}{Theorem}
\newtheorem{lemma}{Lemma}

\newtheorem{remark}[theorem]{Remark}

\newtheorem{example}[theorem]{Example}

\usepackage{lipsum}

\usepackage{amsthm}
\usepackage{comment}
\usepackage[latin1]{inputenc}
\usepackage{subfigure}
\usepackage{color}
\usepackage{amsmath}
\usepackage{amsthm}
\usepackage{amstext}
\usepackage{amssymb}
\usepackage{amsfonts}
\usepackage{graphicx}
\usepackage{young}
\usepackage{multicol}
\usepackage{mathrsfs}
\usepackage[all]{xy}
\usepackage{tikz}
\usepackage{mathtools}
\usepackage{tabularx}
\usepackage{array}
\usepackage{commath}
\usepackage{ytableau}
\usepackage[colorinlistoftodos]{todonotes}

\usepackage{enumerate}

\newcommand{\bd}{\textbf{d}}
\newcommand{\bl}{{\boldsymbol\lambda}}

\newcommand{\fk}{\Bbbk}

\newcommand{\pos}{\textsf{P}_{{Q,m}}}

\def\RSK{\rho_{Q,m}(M)}
\newcommand{\cat}{\mathcal{C}_{Q,m}}
\newcommand{\pro}{\operatorname{pro}}

\newcommand{\JFg}{\operatorname{GenJF}}

\renewcommand{\P}{\sf P}
\newcommand{\J}{\sf J}
\newcommand{\OO}{\sf O}
\newcommand{\x}{\sf x}
\newcommand{\y}{\sf y}

\title[Minuscule Reverse Plane Partitions]{Minuscule Reverse Plane Partitions via Quiver Representations}
\author{Alexander Garver}
\author{Rebecca Patrias}
\author{Hugh Thomas}

\address{LaCIM, Universit\'e du Qu\'ebec \`a Montr\'eal, Qu\'ebec (Canada)}

\begin{document}

\begin{abstract}{The Hillman--Grassl correspondence is a well-known bijection between multisets of rim hooks of a partition shape $\lambda$ and reverse plane partitions of $\lambda$. We use the tools of quiver representations to generalize Hillman--Grassl in type $A$ and to define an analogue in all minuscule types.}\end{abstract}






\maketitle
A reverse plane partition of a Young diagram is a filling of the boxes with $\{0,1,2,\ldots\}$ such that the entries in rows and columns are weakly increasing. They are prominent combinatorial objects with connections to areas like symmetric functions and representation theory (see for example \cite{stanley1971theory}). Their generating function was discovered by Stanley and is as follows, where $h(u)$ denotes the hook length of the cell $u$ in partition $\lambda$ and $|\rho|$ denotes the sum of the entries in reverse plane partition $\rho$.

\begin{theorem}\cite{stanley1971theory}\label{thm:genfunction}
The generating function for reverse plane partitions of shape $\lambda$ with respect to the sum of its entries is 
\[\sum_{\rho}q^{|\rho|}=\prod_{u\in\lambda}\frac{1}{1-q^{h(u)}}.\]
\end{theorem}

The first bijective proof of this generating function was found by Hillman and Grassl in \cite{hillman1976reverse}. The authors give a bijection between nonnegative integer arrays of shape $\lambda$---which can be thought of as multisets of rim hooks of $\lambda$---and reverse plane partitions of $\lambda$, which is now known as the \textit{Hillman--Grassl correspondence}. This correspondence has since been well studied, for example by Gansner in \cite{gansner1981hillman} and by Morales, Pak, and Panova in \cite{morales2015hook}. Another such bijection was given by Pak in \cite{pak2001hook} and rediscovered in a different form by Sulzgruber in \cite{sulzgruberfull}.

More generally, one can define a reverse plane partition on any poset $\textsf{P}$ to be an order-reversing map $\rho:\textsf{P}\to\{0,1,2,\ldots\}$. In the case where $\textsf{P}$ is a minuscule poset, Proctor \cite{proctor1984bruhat} proves the following analogous generating function for reverse plane partitions of $\textsf{P}$ with respect to the sum of its entries: 
\[\sum_{\rho\in RPP(\textsf{P})} q^{|\rho|} =\prod_{\textsf{u}\in \textsf{P}} \frac{1}{1-q^{\text{rank}(u)}}.\] 
In this paper, we use quiver representations to give a uniform bijective proof of this generating function for all minuscule types that generalizes the Hillman--Grassl and Pak correspondences in type $A$. Our bijection is between certain isomorphism classes of representations of a simply-laced Dynkin quiver and reverse plane partitions on the corresponding minuscule poset; see Theorems~\ref{thm:bijection} and \ref{th-rpp}. The interested reader may see \cite{garver2017+} for a full version of this paper with proofs included. 

In Theorem~\ref{thm:bijection}, we describe our bijection as a sequence of piecewise-linear maps that involve an operation known as \textit{toggling}. There is an action on the set of reverse plane partitions known as \textit{promotion} defined using toggles. An important application of our work is that promotion applied $h$ times to  a reverse plane partition of a minuscule poset is the identity where $h$ is the Coxeter number of the associated root system.



\section{Quiver representations}
We briefly review the basics of quiver representations; see \cite{assem2006elements} for more details. A \textit{quiver} $Q$ is a directed graph. In other words, $Q$ is a 4-tuple $(Q_0,Q_1,s,t)$, where $Q_0$ is a set of \textit{vertices}, $Q_1$ is a set of \textit{arrows}, and $s, t:Q_1 \to Q_0$ are two functions defined so that for every $a \in Q_1$, we have $s(a) \xrightarrow{a} t(a)$. In all of our results, we will assume that $Q$ is a \textit{Dynkin quiver}, meaning that its underlying graph is a simply-laced Dynkin diagram.

A \textit{representation} $V = ((V_i)_{i \in Q_0}, (f_a)_{a \in Q_1})$ of a quiver $Q$  is an assignment of a finite-dimensional $\fk$-vector space $V_i$ to each vertex $i$ and a $\fk$-linear map $f_a: V_{s(a)} \rightarrow V_{t(a)}$ to each arrow $a$ where $\fk$ is a field. The \textit{dimension vector} of $V$ is the vector $\textbf{dim}(V):=(\dim V_i)_{i\in Q_0}$. In all of our results, we will assume that $\Bbbk$ is algebraically closed.


\begin{figure}
    \centering
\begin{minipage}{.5\textwidth}
\begin{tikzpicture}
\node (1) at (0,0) {0};
\node (2) at (2,0) {0};
\node (3) at (4,0) {$\Bbbk$};
\node (4) at (6,0) {0};
\draw[->] (1)--(2) node [midway, above, fill=white] {\small{0}};
\draw[->] (2)--(3) node [midway, above, fill=white] {\small{0}};
\draw[->] (4)--(3) node [midway, above, fill=white] {\small{0}};
\end{tikzpicture}\hspace{.3in}\raisebox{.05in}{0010}\\
\begin{tikzpicture}
\node (1) at (0,0) {0};
\node (2) at (2,0) {0};
\node (3) at (4,0) {$\Bbbk$};
\node (4) at (6,0) {$\Bbbk$};
\draw[->] (1)--(2) node [midway, above, fill=white] {\small{0}};
\draw[->] (2)--(3) node [midway, above, fill=white] {\small{0}};
\draw[->] (4)--(3) node [midway, above, fill=white] {\small{1}};
\end{tikzpicture}\hspace{.3in}\raisebox{.05in}{0011}\\
\begin{tikzpicture}
\node (1) at (0,0) {0};
\node (2) at (2,0) {$\Bbbk$};
\node (3) at (4,0) {$\Bbbk$};
\node (4) at (6,0) {0};
\draw[->] (1)--(2) node [midway, above, fill=white] {\small{0}};
\draw[->] (2)--(3) node [midway, above, fill=white] {\small{1}};
\draw[->] (4)--(3) node [midway, above, fill=white] {\small{0}};
\end{tikzpicture}\hspace{.3in}\raisebox{.05in}{0110}\\
\begin{tikzpicture}
\node (1) at (0,0) {$\Bbbk$};
\node (2) at (2,0) {$\Bbbk$};
\node (3) at (4,0) {$\Bbbk$};
\node (4) at (6,0) {0};
\draw[->] (1)--(2) node [midway, above, fill=white] {\small{1}};
\draw[->] (2)--(3) node [midway, above, fill=white] {\small{1}};
\draw[->] (4)--(3) node [midway, above, fill=white] {\small{0}};
\end{tikzpicture}\hspace{.3in}\raisebox{.05in}{1110}\\
\begin{tikzpicture}
\node (1) at (0,0) {0};
\node (2) at (2,0) {$\Bbbk$};
\node (3) at (4,0) {$\Bbbk$};
\node (4) at (6,0) {$\Bbbk$};
\draw[->] (1)--(2) node [midway, above, fill=white] {\small{0}};
\draw[->] (2)--(3) node [midway, above, fill=white] {\small{1}};
\draw[->] (4)--(3) node [midway, above, fill=white] {\small{1}};
\end{tikzpicture}\hspace{.3in}\raisebox{.05in}{0111}\\
\begin{tikzpicture}
\node (1) at (0,0) {$\Bbbk$};
\node (2) at (2,0) {$\Bbbk$};
\node (3) at (4,0) {$\Bbbk$};
\node (4) at (6,0) {$\Bbbk$};
\draw[->] (1)--(2) node [midway, above, fill=white] {\small{1}};
\draw[->] (2)--(3) node [midway, above, fill=white] {\small{1}};
\draw[->] (4)--(3) node [midway, above, fill=white] {\small{1}};
\end{tikzpicture}\hspace{.3in}\raisebox{.05in}{1111}
\end{minipage}
\begin{minipage}{.4\textwidth}
\scalebox{.8}{
   \begin{tikzpicture}
\node(0011) at (0,0) {0011};
\node(0010) at (-1,1) {0010};
\node(0111) at (1,1) {0111};
\node(0110) at (0,2) {0110};
\node(1111) at (2,2) {1111};
\node(1110) at (1,3) {1110};
\node(0100) at (2,0) {\textcolor{gray}{0100}};
\node(1000) at (4,0) {\textcolor{gray}{1000}};
\node(1100) at (3,1) {\textcolor{gray}{1100}};
\node(0001) at (3,3) {\textcolor{gray}{0001}};
\draw[->] (0010)--(0011);
\draw[->] (0010) -- (0110);
\draw[->] (0011) -- (0111);
\draw[->] (0110) -- (0111);
\draw[->] (0110) -- (1110);
\draw[->] (0111) -- (1111);
\draw[->] (1110) -- (1111);
\draw[gray,dashed,->] (0111)--(0100);
\draw[gray,dashed,->] (0100)--(1100);
\draw[gray,dashed,->] (1100)--(1000);
\draw[gray,dashed,->] (1111)--(1100);
\draw[gray,dashed,->] (1111)--(0001);
\draw[red,dotted,->] (1000)--(0100) node [midway, above, fill=white] {\small{$\tau$}};
\draw[red,dotted,->] (0100)--(0011) node [midway, above, fill=white] {\small{$\tau$}};
\draw[red,dotted,->] (1100)--(0111) node [midway, above, fill=white] {\small{$\tau$}};
\draw[red,dotted,->] (0111)--(0010) node [midway, above, fill=white] {\small{$\tau$}};
\draw[red,dotted,->] (1111)--(0110) node [midway, above, fill=white] {\small{$\tau$}};
\draw[red,dotted,->] (0001)--(1110) node [midway, above, fill=white] {\small{$\tau$}};
\end{tikzpicture}}
\end{minipage}
    \caption{Let $Q=1\rightarrow2\rightarrow3\leftarrow4$ with $m=3$. On the left is the set of indecomposable representations (up to isomorphism) of $Q$ that have support over vertex 3 and their dimension vectors. On the right is the minuscule poset (in black) inside of the AR quiver of $Q$. The action of $\tau$ is shown in red; the leftmost representations map to 0.}
    \label{fig:irreps}
\end{figure}
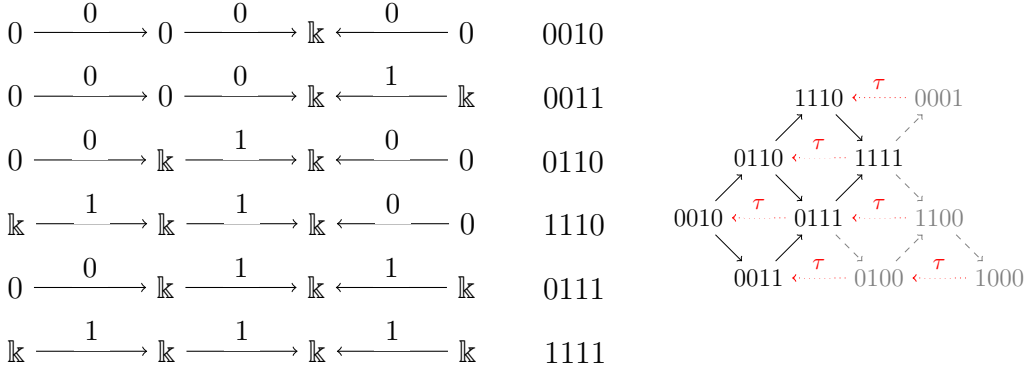

Let $V = ((V_i)_{i \in Q_0}, (f_a)_{a \in Q_1})$ and $W  = ((W_i)_{i \in Q_0}, (g_a)_{a \in Q_1})$ be two representations of a quiver $Q$. 
A \textit{morphism} $\theta : V \rightarrow W$ consists of a collection of linear maps $\theta_i : V_i \rightarrow W_i$ that are compatible with each of the linear maps in $V$ and $W$.  That is, for each arrow $a \in Q_1$, we have $\theta_{t(a)} \circ f_a = g_a \circ \theta_{s(a)}$. We let $\text{Hom}(V,W)$ denote the $\Bbbk$-vector space of all morphisms $\theta: V \to W.$ Additionally, we say that a collection of linear maps $\{\theta_i: V_i \to V_i\}_{i\in Q_0}$ is \textit{compatible} with a representation $V$ if they define a morphism of $V$. 

The representations of a quiver $Q$ along with morphisms between them form an abelian category, which is equivalent to the category of finitely generated left modules over the path algebra of $Q$. It is therefore natural to ask for a classification of the indecomposable representations of $Q$, i.e., those which cannot be written as a direct sum of two non-zero representations. Since we will restrict to the case when $Q$ is Dynkin, Gabriel's Theorem \cite{gabriel1972unzerlegbare} implies that the indecomposable representations of $Q$ are in one-to-one correspondence with the positive roots of the root system corresponding to $Q$, under the map sending a representation to its dimension vector.  Figure~\ref{fig:irreps} shows examples of indecomposable representations of a Dynkin quiver of type $A_4$ with their dimension vectors.




One can understand a great deal about the category of representations of $Q$ by constructing the \textit{Auslander--Reiten quiver} (or \textit{AR quiver}) of $Q$. In this paper, we use the AR quiver of $Q$ to relate the representation theory of $Q$ to the combinatorics of minuscule posets. By definition, the \textit{AR quiver} of $Q$ is the quiver whose vertices are indexed by the isomorphism classes of indecomposable representations of $Q$ and whose arrows index a basis of the space of irreducible morphisms between the corresponding representations. In Figure~\ref{fig:irreps}, we show an example of an AR quiver. 

There is a well-known endofunctor on the category of representations of $Q$, called the \textit{Auslander--Reiten translation} and denoted by $\tau$. We omit its homological definition, but we show how $\tau$ acts on indecomposables in Figure~\ref{fig:irreps}. For our work, the important fact is that $\tau$ partitions the indecomposable representations of $Q$ into  $\tau$\textit{-orbits}.


\section{Minuscule posets}
\begin{figure}
\begin{minipage}{.5\textwidth}
\begin{displaymath}
    \scalebox{.75}{\xymatrix{A_n & \textcolor{red}{1} \ar@{-}[r] & \textcolor{red}{2} \ar@{-}[r] & & \cdots & & \textcolor{red}{n}\ar@{-}[l] & & \\
    & & & & & & \textcolor{red}{n} \ar@{-}[dl] \\
    D_n & \textcolor{red}{1} \ar@{-}[r] & 2 \ar@{-}[r] & & \cdots & n-2 \\
    & & & 6 \ar@{-}[d] & & & \textcolor{red}{n-1}\ar@{-}[ul] \\
    E_6 & \textcolor{red}{1}\ar@{-}[r] & 2 \ar@{-}[r] & 3 \ar@{-}[r] & 4 \ar@{-}[r] & \textcolor{red}{5} \\
   & & & 7 \ar@{-}[d] \\
    E_7 & 1\ar@{-}[r] & 2 \ar@{-}[r] & 3 \ar@{-}[r] & 4 \ar@{-}[r] & 5 \ar@{-}[r] & \textcolor{red}{6} \\}}
\end{displaymath}
\end{minipage}
\begin{minipage}{.5\textwidth}
$$\begin{array}{ccccccccccc} \includegraphics[scale=.55]{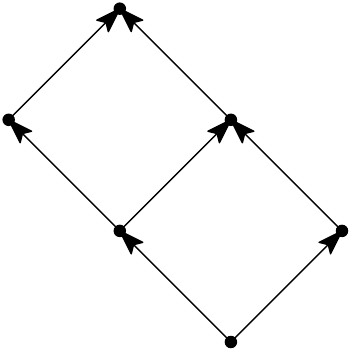} & &  \includegraphics[scale=.55]{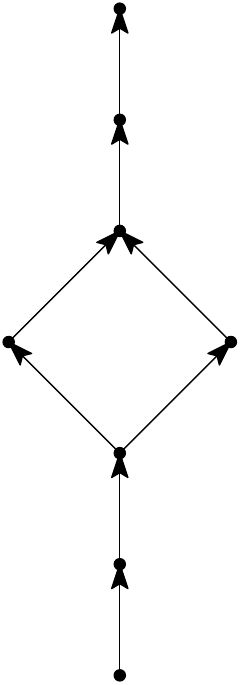}  & & \includegraphics[scale=.55]{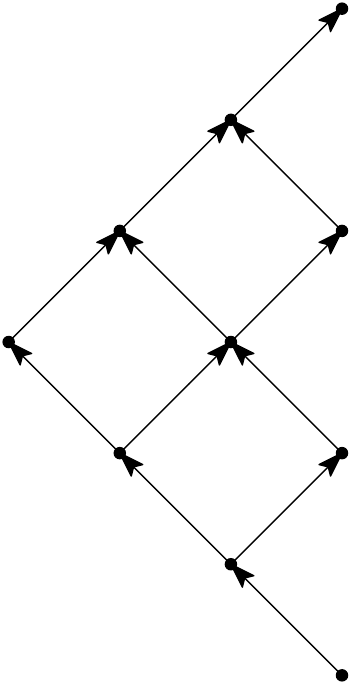}\\ (a)  & &  (b)  & &  (c)   \end{array}$$
\end{minipage}
\caption{The possible choices of minuscule vertices of Dynkin diagrams are shown on the left in red. On the right, we have the minuscule poset of  type $A_4$ with $m = 3$, of type $D_5$ with $m=1$, and type $D_5$ with $m = 4$, respectively, in ($a$), ($b$), and ($c$).  }
\label{diagrams}
\end{figure}

A vertex $m$ of $Q$ is called 
\emph{minuscule} if every positive root in the root system of $Q$ has the simple root corresponding to vertex $m$ appearing in its simple root expansion with coefficient 0 or 1.
In other words, a vertex $m$ of $Q$ is minuscule if every indecomposable representation of $Q$ supported over $m$ has dimension 1 at $m$.  For example, every vertex of a type $A$ Dynkin quiver is minuscule. Figure~\ref{diagrams} shows the possible choices of minuscule vertices.





Each choice of minuscule vertex gives rise to a minuscule poset. We now give explicit descriptions of the minuscule posets, which appear in the context of minuscule representations of complex semisimple Lie algebras and were classified up to isomorphism by Proctor in \cite{proctor1984bruhat}. Recall that for a poset $\P$, we let $\J(\P)$ denote its poset of \textit{order ideals} (i.e., subsets of $\OO \subset \P$ where if $\x \in \OO$ and $\y \le_{\P} \x$, one has that $\y \in \OO$), ordered by inclusion.  By \cite[Theorem 8.3.10]{green2013combinatorics}, there is a minuscule poset for each choice of a simply-laced Dynkin diagram other than $E_8$ with a chosen minuscule vertex. Their isomorphism types appear in Table~\ref{table:1}. 
There, we write $[n]$ for the poset that is a chain whose elements are $1, \ldots, n$ in increasing order. Figure~\ref{diagrams} shows some examples of minuscule posets.

\begin{table}[!htbp]
\centering
\begin{tabular}{c c c } 
 \hline
 Type & $m$ & minuscule poset \\ [0.5ex] 
 \hline\hline
 $A_n$ & $k$ & $[k]\times[n+1-k]$ \\ 
 $D_n$ & $1$ & $\J^{n-3}([2]\times[2])$ \\
 $D_n$ & $n-1$, $n$ & $\J([2]\times[n-2])$ \\
 $E_6$ & $1$, $5$ & $\J^2([2]\times[3])$ \\
 $E_7$ & $6$ & $\J^3([2]\times[3])$  \\ [1ex] 
 \hline
\end{tabular}
\caption{The isomorphism types of the minuscule posets. Here we are referring to the vertex labeling of the Dynkin diagrams appearing in Figure~\ref{diagrams}.}
\label{table:1}
\end{table}

Our next result shows how the minuscule posets are related to the representation theory of quivers. Given a minuscule vertex $m \in Q_0$, let $\pos$ denote the poset that is the transitive closure of the arrows in the full subquiver of the AR quiver of $Q$ whose vertex set is the set of isomorphism classes of representations supported at $m$. Examples for type $A$ quivers are shown in Figures~\ref{fig:irreps}, \ref{fig:HGexample}, and \ref{fig:Pakexample}. We let $\cat$ denote the category consisting of all direct sums of indecomposable representations of $Q$ supported at $m$.

{We will need the following lemma, which is not hard to verify but appears to be new.} 

\begin{lemma}
The poset $\pos$ is isomorphic to the minuscule poset determined by the underlying graph of $Q$ and minuscule vertex $m$.
\end{lemma}

\section{Quiver representations and reverse plane partitions}
Our main result is the following theorem. 

\begin{theorem}\label{thm:bijection}
Let $Q$ be a Dynkin quiver with minuscule vertex $m$. There is a bijection between isomorphism classes of objects in $\mathcal{C}_{Q,m}$ and reverse plane partitions of the minuscule poset $\pos$.
\end{theorem}

\begin{remark}
It is also possible to give a bijection between reverse plane partitions in any order filter of $\pos$ and the isomorphism classes of objects in a suitable category.  For simplicity, we only discuss the bijection for the full posets here.
\end{remark}

\begin{remark}
{Theorem~\ref{thm:bijection} provides a proof of Proctor's generating function identity from the Introduction. We may re-express that equation as follows, where the first equality is implied by Theorem~\ref{thm:bijection} and the second holds since any representation in $\mathcal{C}_{Q,m}$ uniquely decomposes as a direct sum of indecomposable representations $M_{\textsf{u}}$ belonging to $\mathcal{C}_{Q,m}$:
\[\sum_{\rho\in RPP(\textsf{P}_{Q,m})} q^{|\rho|} = \sum_{V\in \mathcal{C}_{Q,m}} q^{\dim V} =\prod_{\textsf{u}\in \textsf{P}_{Q,m}} \frac{1}{1-q^{\dim M_{\textsf{u}}}}.\] In \cite[Corollary 4.9]{garver2017+}, we show that there is a multivariate version of this generating function identity that recovers an identity of Gansner \cite[Corollary 5.2]{gansner1981hillman}.}
\end{remark}

In Section~\ref{sec:combdesc}, we describe the bijection explicitly in terms of \textit{toggling}, which we introduce in the next section. In Section~\ref{sec:conceptualdescription}, we give an alternative description of the bijection using the geometry of representation varieties.

\subsection{Toggling}
Let $\rho: \textsf{P} \to [0,N]$ be a reverse plane partition.  
For $\textsf{x}$ an element of $\textsf{P}$, we define the \textit{toggle} of $\rho$ at $\textsf{x} \in \textsf{P}$ by
$$\begin{array}{rcl}
t_\textsf{x}\rho(\textsf{y}) & = & \left\{\begin{array}{lcl} \displaystyle \max_{\textsf{y} \lessdot \textsf{y}_1}\rho(\textsf{y}_1) +\min_{\textsf{y}_2\lessdot \textsf{y}}\rho(\textsf{y}_2)-\rho(\textsf{y}) & : & \text{if $\textsf{y} = \textsf{x}$} \\ \rho(\textsf{y}) & : & \text{if $\textsf{y} \neq \textsf{x}$,}\end{array}\right.
\end{array}$$
where $\textsf{y}$ is any element of $\textsf{P}$.  If $\textsf{y}$ is maximal,
we interpret
$\max_{\textsf{y} \lessdot \textsf{y}_1}\rho(\textsf{y}_1)$ as 0, and if $\textsf{y}$ is minimal, we 
interpret $\min_{\textsf{y}_2\lessdot \textsf{y}}\rho(\textsf{y}_2)$ as $N$.
Since $\rho$ is a reverse plane partition, so is $t_\textsf{x} \rho.$ Additionally, observe that $t_\textsf{x}\circ t_\textsf{x}(\rho) = \rho$. This is the toggle operation considered by Einstein and Propp (up to rescaling, and restricted to lattice points) \cite{einstein2013combinatorial}.

\subsection{Combinatorial Description}\label{sec:combdesc}


We now give an explicit description of the bijection in Theorem~\ref{thm:bijection}, which we will denote by $\rho_{Q,m}(-)$. We must first establish a linear order on the indecomposable representations of $Q$. Choose a linear order on the indecomposable representations of $Q$ compatible with the opposite of the AR quiver order. Let $M_1,\dots, M_s$ denote the indecomposable representations of $Q$ in this order.



Let $M=\bigoplus_{j=1}^s M_j^{c_j} \in \cat.$  The reverse plane partition $\RSK$ corresponding to $M$ is obtained by constructing a sequence of fillings of the minuscule poset $\pos$, starting with $\rho_0$, which is the zero filling. These fillings are defined by
$$\begin{array}{ccccc}\rho_k(\textsf{x}) & := & \left\{\begin{array}{lcl} \displaystyle \max_{\textsf{x}\lessdot \textsf{y}}\rho_{k-1}(\textsf{y}) + c_k & : & \text{if $\textsf{x}$ is the element of 
$\pos$ corresponding to $M_k$,} \\
({t}_\textsf{x}\rho_{k-1})(\textsf{x}) & : & \text{if $\textsf{x}$ corresponds to  $\tau^\ell(M_k)$ for some $\ell < 0$, and}\\
\rho_{k-1}(\textsf{x}) & : & \text{otherwise,}
\end{array}\right.\end{array}$$ 
where $\textsf{x}$ is any element of $\pos$. Then we have that $\RSK = \rho_s.$ 


Note that using this description of the algorithm, the intermediate fillings $\rho_k$ for $k<s$ are not reverse plane partitions of $\pos$. However, by restricting $\rho_k$ to the elements of $\pos$ corresponding to $M_1,\ldots,M_k$, we do obtain a reverse plane partition on the induced subposet of $\pos$ whose elements correspond to $M_1,\ldots,M_k$.

\begin{figure}
\centering
\scalebox{.8}{
\begin{tikzpicture}[scale=1]
\node (11100) at (1,0) {$00100_{6}$};
\node (00010) at (3,0) {\textcolor{gray}{$00010_3$}};
\node (00001) at (5,0) {\textcolor{gray}{$00001_1$}}; 
\node (01100) at (0,1) {$01100_{9}$};
\node (11110) at (2,1) {$00110_5$};
\node (00011) at (4,1) {\textcolor{gray}{$00011_2$}};
\node (01110) at (1,2) {$01110_{8}$};
\node (11111) at (3,2) {$00111_4$};
\node (01111) at (2,3) {$01111_7$};
\node (00100) at (-1,2) {$11100_{13}$};
\node (00110) at (0,3) {$11110_{12}$};
\node (00111) at (1,4) {$11111_{11}$};
\node (10000) at (-3,0) {\textcolor{gray}{$10000_{15}$}};
\node (01000) at (-1,0) {\textcolor{gray}{$01000_{10}$}};
\node (11000) at (-2,1) {\textcolor{gray}{$11000_{14}$}};
\draw[->] (00100)--(01100);
\draw[->] (01100)--(11100);
\draw[->] (00100)--(00110);
\draw[->] (00110)--(00111);
\draw[->] (00110)--(01110);
\draw[->] (01100)--(01110);
\draw[->] (11100)--(11110);
\draw[->] (11110)--(11111);
\draw[->] (10000)--(11000);
\draw[->] (01110)--(11110);
\draw[->] (11110)--(00010);
\draw[->] (00111)--(01111);
\draw[->] (01111)--(11111);
\draw[->] (11111)--(00011);
\draw[->] (00011)--(00001);
\draw[->] (00010)--(00011);
\draw[->] (01110)--(01111);
\draw[->] (11000)--(01000);
\draw[->] (11000)--(00100);
\draw[->] (01000)--(01100);
\end{tikzpicture}
\ytableausetup{boxsize=.2in}
\raisebox{1in}{$M=$
\raisebox{-.65in}{\begin{tikzpicture}[scale=.7]
\node (11100) at (1,0) {0};
\node (01100) at (0,1) {3};
\node (11110) at (2,1) {1};
\node (01110) at (1,2) {1};
\node (11111) at (3,2) {1};
\node (01111) at (2,3) {0};
\node (00100) at (-1,2) {4};
\node (00110) at (0,3) {0};
\node (00111) at (1,4) {2};
\draw[->] (00100)--(01100);
\draw[->] (01100)--(11100);
\draw[->] (00100)--(00110);
\draw[->] (00110)--(00111);
\draw[->] (00110)--(01110);
\draw[->] (01100)--(01110);
\draw[->] (11100)--(11110);
\draw[->] (11110)--(11111);
\draw[->] (01110)--(11110);
\draw[->] (00111)--(01111);
\draw[->] (01111)--(11111);
\draw[->] (01110)--(01111);
\end{tikzpicture}}}}
\scalebox{.8}{
\begin{tikzpicture}[scale=.7]
\node (11100) at (1,0) {\textcolor{lightgray}{0}};
\node (01100) at (0,1) {\textcolor{lightgray}{0}};
\node (11110) at (2,1) {\textcolor{lightgray}{0}};
\node (01110) at (1,2) {\textcolor{lightgray}{0}};
\node (11111) at (3,2) {1};
\node (01111) at (2,3) {\textcolor{lightgray}{0}};
\node (00100) at (-1,2) {\textcolor{lightgray}{0}};
\node (00110) at (0,3) {\textcolor{lightgray}{0}};
\node (00111) at (1,4) {\textcolor{lightgray}{0}};
\draw[->] (00100)--(01100);
\draw[->] (01100)--(11100);
\draw[->] (00100)--(00110);
\draw[->] (00110)--(00111);
\draw[->] (00110)--(01110);
\draw[->] (01100)--(01110);
\draw[->] (11100)--(11110);
\draw[->] (11110)--(11111);
\draw[->] (01110)--(11110);
\draw[->] (00111)--(01111);
\draw[->] (01111)--(11111);
\draw[->] (01110)--(01111);
\end{tikzpicture}
\begin{tikzpicture}[scale=.7]
\node (11100) at (1,0) {\textcolor{lightgray}{0}};
\node (01100) at (0,1) {\textcolor{lightgray}{0}};
\node (11110) at (2,1) {2};
\node (01110) at (1,2) {\textcolor{lightgray}{0}};
\node (11111) at (3,2) {1};
\node (01111) at (2,3) {\textcolor{lightgray}{0}};
\node (00100) at (-1,2) {\textcolor{lightgray}{0}};
\node (00110) at (0,3) {\textcolor{lightgray}{0}};
\node (00111) at (1,4) {\textcolor{lightgray}{0}};
\draw[->] (00100)--(01100);
\draw[->] (01100)--(11100);
\draw[->] (00100)--(00110);
\draw[->] (00110)--(00111);
\draw[->] (00110)--(01110);
\draw[->] (01100)--(01110);
\draw[->] (11100)--(11110);
\draw[->] (11110)--(11111);
\draw[->] (01110)--(11110);
\draw[->] (00111)--(01111);
\draw[->] (01111)--(11111);
\draw[->] (01110)--(01111);
\end{tikzpicture}
\begin{tikzpicture}[scale=.7]
\node (11100) at (1,0) {2};
\node (01100) at (0,1) {\textcolor{lightgray}{0}};
\node (11110) at (2,1) {2};
\node (01110) at (1,2) {\textcolor{lightgray}{0}};
\node (11111) at (3,2) {1};
\node (01111) at (2,3) {\textcolor{lightgray}{0}};
\node (00100) at (-1,2) {\textcolor{lightgray}{0}};
\node (00110) at (0,3) {\textcolor{lightgray}{0}};
\node (00111) at (1,4) {\textcolor{lightgray}{0}};
\draw[->] (00100)--(01100);
\draw[->] (01100)--(11100);
\draw[->] (00100)--(00110);
\draw[->] (00110)--(00111);
\draw[->] (00110)--(01110);
\draw[->] (01100)--(01110);
\draw[->] (11100)--(11110);
\draw[->] (11110)--(11111);
\draw[->] (01110)--(11110);
\draw[->] (00111)--(01111);
\draw[->] (01111)--(11111);
\draw[->] (01110)--(01111);
\end{tikzpicture}
\begin{tikzpicture}[scale=.7]
\node (11100) at (1,0) {2};
\node (01100) at (0,1) {\textcolor{lightgray}{0}};
\node (11110) at (2,1) {2};
\node (01110) at (1,2) {\textcolor{lightgray}{0}};
\node (11111) at (3,2) {1};
\node (01111) at (2,3) {1};
\node (00100) at (-1,2) {\textcolor{lightgray}{0}};
\node (00110) at (0,3) {\textcolor{lightgray}{0}};
\node (00111) at (1,4) {\textcolor{lightgray}{0}};
\draw[->] (00100)--(01100);
\draw[->] (01100)--(11100);
\draw[->] (00100)--(00110);
\draw[->] (00110)--(00111);
\draw[->] (00110)--(01110);
\draw[->] (01100)--(01110);
\draw[->] (11100)--(11110);
\draw[->] (11110)--(11111);
\draw[->] (01110)--(11110);
\draw[->] (00111)--(01111);
\draw[->] (01111)--(11111);
\draw[->] (01110)--(01111);
\end{tikzpicture}
\begin{tikzpicture}[scale=.7]
\node (11100) at (1,0) {2};
\node (01100) at (0,1) {\textcolor{lightgray}{0}};
\node (11110) at (2,1) {2};
\node (01110) at (1,2) {3};
\node (11111) at (3,2) {0};
\node (01111) at (2,3) {1};
\node (00100) at (-1,2) {\textcolor{lightgray}{0}};
\node (00110) at (0,3) {\textcolor{lightgray}{0}};
\node (00111) at (1,4) {\textcolor{lightgray}{0}};
\draw[->] (00100)--(01100);
\draw[->] (01100)--(11100);
\draw[->] (00100)--(00110);
\draw[->] (00110)--(00111);
\draw[->] (00110)--(01110);
\draw[->] (01100)--(01110);
\draw[->] (11100)--(11110);
\draw[->] (11110)--(11111);
\draw[->] (01110)--(11110);
\draw[->] (00111)--(01111);
\draw[->] (01111)--(11111);
\draw[->] (01110)--(01111);
\end{tikzpicture}}
\scalebox{.8}{
\begin{tikzpicture}[scale=.7]
\node (11100) at (1,0) {2};
\node (01100) at (0,1) {6};
\node (11110) at (2,1) {0};
\node (01110) at (1,2) {3};
\node (11111) at (3,2) {0};
\node (01111) at (2,3) {1};
\node (00100) at (-1,2) {\textcolor{lightgray}{0}};
\node (00110) at (0,3) {\textcolor{lightgray}{0}};
\node (00111) at (1,4) {\textcolor{lightgray}{0}};
\draw[->] (00100)--(01100);
\draw[->] (01100)--(11100);
\draw[->] (00100)--(00110);
\draw[->] (00110)--(00111);
\draw[->] (00110)--(01110);
\draw[->] (01100)--(01110);
\draw[->] (11100)--(11110);
\draw[->] (11110)--(11111);
\draw[->] (01110)--(11110);
\draw[->] (00111)--(01111);
\draw[->] (01111)--(11111);
\draw[->] (01110)--(01111);
\end{tikzpicture}
\begin{tikzpicture}[scale=.7]
\node (11100) at (1,0) {4};
\node (01100) at (0,1) {6};
\node (11110) at (2,1) {0};
\node (01110) at (1,2) {3};
\node (11111) at (3,2) {0};
\node (01111) at (2,3) {1};
\node (00100) at (-1,2) {\textcolor{lightgray}{0}};
\node (00110) at (0,3) {\textcolor{lightgray}{0}};
\node (00111) at (1,4) {\textcolor{lightgray}{0}};
\draw[->] (00100)--(01100);
\draw[->] (01100)--(11100);
\draw[->] (00100)--(00110);
\draw[->] (00110)--(00111);
\draw[->] (00110)--(01110);
\draw[->] (01100)--(01110);
\draw[->] (11100)--(11110);
\draw[->] (11110)--(11111);
\draw[->] (01110)--(11110);
\draw[->] (00111)--(01111);
\draw[->] (01111)--(11111);
\draw[->] (01110)--(01111);
\end{tikzpicture}
\begin{tikzpicture}[scale=.7]
\node (11100) at (1,0) {4};
\node (01100) at (0,1) {6};
\node (11110) at (2,1) {0};
\node (01110) at (1,2) {3};
\node (11111) at (3,2) {0};
\node (01111) at (2,3) {1};
\node (00100) at (-1,2) {\textcolor{lightgray}{0}};
\node (00110) at (0,3) {\textcolor{lightgray}{0}};
\node (00111) at (1,4) {3};
\draw[->] (00100)--(01100);
\draw[->] (01100)--(11100);
\draw[->] (00100)--(00110);
\draw[->] (00110)--(00111);
\draw[->] (00110)--(01110);
\draw[->] (01100)--(01110);
\draw[->] (11100)--(11110);
\draw[->] (11110)--(11111);
\draw[->] (01110)--(11110);
\draw[->] (00111)--(01111);
\draw[->] (01111)--(11111);
\draw[->] (01110)--(01111);
\end{tikzpicture}
\begin{tikzpicture}[scale=.7]
\node (11100) at (1,0) {4};
\node (01100) at (0,1) {6};
\node (11110) at (2,1) {0};
\node (01110) at (1,2) {3};
\node (11111) at (3,2) {0};
\node (01111) at (2,3) {2};
\node (00100) at (-1,2) {\textcolor{lightgray}{0}};
\node (00110) at (0,3) {3};
\node (00111) at (1,4) {3};
\draw[->] (00100)--(01100);
\draw[->] (01100)--(11100);
\draw[->] (00100)--(00110);
\draw[->] (00110)--(00111);
\draw[->] (00110)--(01110);
\draw[->] (01100)--(01110);
\draw[->] (11100)--(11110);
\draw[->] (11110)--(11111);
\draw[->] (01110)--(11110);
\draw[->] (00111)--(01111);
\draw[->] (01111)--(11111);
\draw[->] (01110)--(01111);
\end{tikzpicture}
\begin{tikzpicture}[scale=.7]
\node (11100) at (1,0) {4};
\node (01100) at (0,1) {6};
\node (11110) at (2,1) {0};
\node (01110) at (1,2) {2};
\node (11111) at (3,2) {0};
\node (01111) at (2,3) {2};
\node (00100) at (-1,2) {10};
\node (00110) at (0,3) {3};
\node (00111) at (1,4) {3};
\draw[->] (00100)--(01100);
\draw[->] (01100)--(11100);
\draw[->] (00100)--(00110);
\draw[->] (00110)--(00111);
\draw[->] (00110)--(01110);
\draw[->] (01100)--(01110);
\draw[->] (11100)--(11110);
\draw[->] (11110)--(11111);
\draw[->] (01110)--(11110);
\draw[->] (00111)--(01111);
\draw[->] (01111)--(11111);
\draw[->] (01110)--(01111);
\end{tikzpicture}}
\scalebox{.8}{
\begin{tikzpicture}[scale=.7]
\node (11100) at (1,0) {4};
\node (01100) at (0,1) {8};
\node (11110) at (2,1) {2};
\node (01110) at (1,2) {2};
\node (11111) at (3,2) {0};
\node (01111) at (2,3) {2};
\node (00100) at (-1,2) {10};
\node (00110) at (0,3) {3};
\node (00111) at (1,4) {3};
\draw[->] (00100)--(01100);
\draw[->] (01100)--(11100);
\draw[->] (00100)--(00110);
\draw[->] (00110)--(00111);
\draw[->] (00110)--(01110);
\draw[->] (01100)--(01110);
\draw[->] (11100)--(11110);
\draw[->] (11110)--(11111);
\draw[->] (01110)--(11110);
\draw[->] (00111)--(01111);
\draw[->] (01111)--(11111);
\draw[->] (01110)--(01111);
\end{tikzpicture}
\begin{tikzpicture}[scale=.7]
\node (11100) at (1,0) {6};
\node (01100) at (0,1) {8};
\node (11110) at (2,1) {2};
\node (01110) at (1,2) {2};
\node (11111) at (3,2) {0};
\node (01111) at (2,3) {2};
\node (00100) at (-1,2) {10};
\node (00110) at (0,3) {3};
\node (00111) at (1,4) {3};
\draw[->] (00100)--(01100);
\draw[->] (01100)--(11100);
\draw[->] (00100)--(00110);
\draw[->] (00110)--(00111);
\draw[->] (00110)--(01110);
\draw[->] (01100)--(01110);
\draw[->] (11100)--(11110);
\draw[->] (11110)--(11111);
\draw[->] (01110)--(11110);
\draw[->] (00111)--(01111);
\draw[->] (01111)--(11111);
\draw[->] (01110)--(01111);
\end{tikzpicture}
\raisebox{.6in}{$=\RSK$}\hspace{.7in}
\raisebox{.8in}{
\begin{ytableau}
0 & 2 & 3 \\
2 & 2 & 3 \\
6 & 8 & 10
\end{ytableau}}
}
\caption[mycaption]{The top left image is the AR quiver associated with the quiver $Q=1\leftarrow 2 \leftarrow \mathbf{3}\leftarrow 4 \leftarrow 5$ with chosen minuscule vertex 3. The dimension vectors with support on vertex 3 are black, while the others are gray. The top right image represents a representation $M$ whose indecomposable summands all have support on vertex 3. The images below show the steps in computing $\RSK$, 
starting at $\rho_4$, 
and the bottom right shows the corresponding reverse plane partition of a Young diagram. Note that these posets are oriented from left to right.
{The computation of $\rho_{8}$ is the first time in the process that we apply a toggle. The computation of $\rho_{10}$ is the first time in the process that we apply only a toggle.}}
\label{fig:HGexample}
\end{figure}
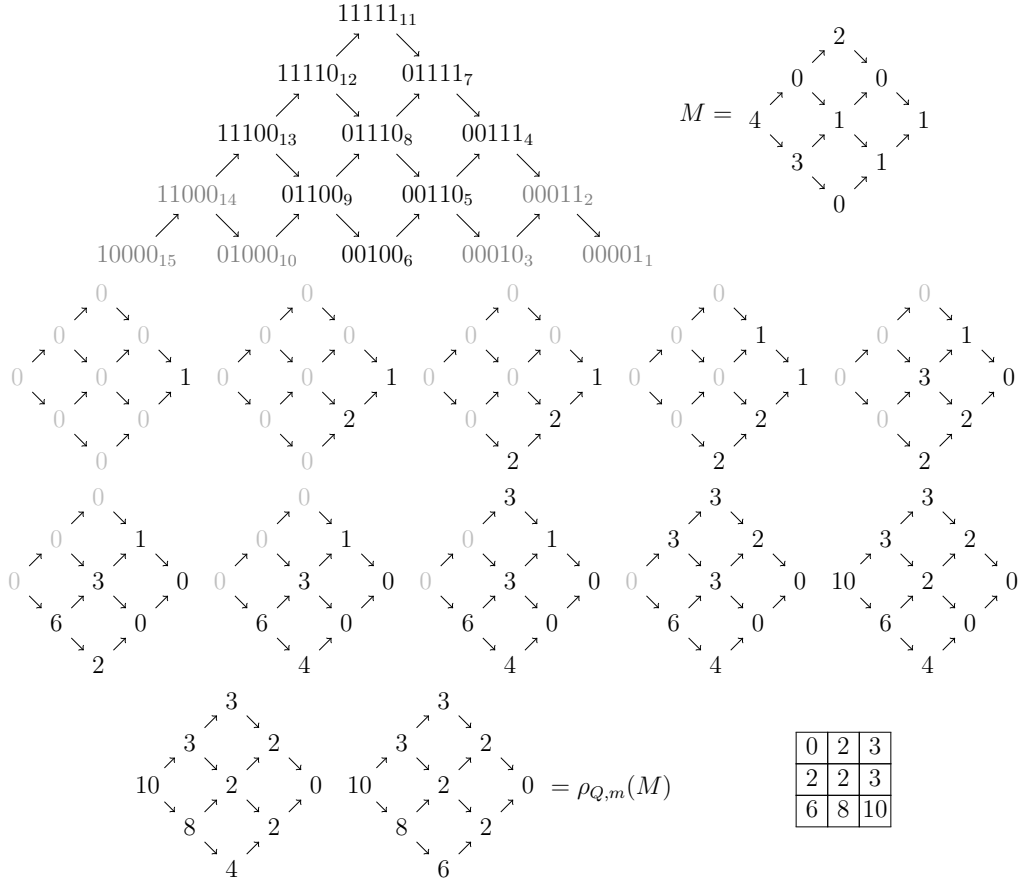

Observe that in the process of constructing $\rho_s$, we never toggle at a minimal element of one of these induced posets. Therefore, the result of the procedure does not depend on the choice of $N$ used in the definition of toggling.
See Figures~\ref{fig:HGexample} and \ref{fig:Pakexample} for examples in type $A$ worked out step-by-step using this explicit description. 

\subsection{Nilpotent endomorphisms of quiver representations}\label{sec:conceptualdescription}
Let $V$ be a 
representation of $Q$ with dimension vector $\bd$, over an algebraically closed ground field $\fk$.  Let $\phi$ be a nilpotent endomorphism of $V$.  At each
vertex $i$ of $Q$, the endomorphism $\phi$ induces an endomorphism $\phi_i$ of $V_i$.  We can consider the
Jordan form of each of these vector space endomorphisms,
which gives us a sequence of partitions $\lambda^i \vdash d_i$.  We show that for a generic choice of $\phi$, the 
$n$-tuple $\bl=(\lambda^1,\dots,\lambda^n)$ is well-defined.
We refer to this as the \textit{Jordan data} of $V$, and we write it as
$\JFg(V)$. The following is an alternate description of the bijection $\rho_{Q,m}(-)$.

\begin{theorem}  \label{th-rpp}
For $V\in \cat$, define a map $\rho_V:\pos \rightarrow \mathbb{N}$ as follows: The values of $\rho_V$ restricted to the $\tau$-orbit corresponding to vertex $j$ are the entries 
of $\JFg(V)^j$, padded with extra zeros if necessary, and ordered so that, restricted to this $\tau$-orbit, the function is order-reversing.
Then $\rho_V$ is a reverse plane partition of $\pos$.
The map from isomorphism classes in $\cat$ to 
reverse plane partitions of $\pos$, sending
$V$ to $\rho_V$, is a bijection that agrees with $\rho_{Q,m}(-).$

\end{theorem}

\begin{example}\label{ex:intro1}
Let $Q=1\rightarrow 2 \leftarrow 3$ and $m=2$. 
By identifying indecomposable representations of $Q$ with their dimension vectors, each $V\in \mathcal{C}_{Q,2}$ is isomorphic to $010^a\oplus 011^b\oplus 110^c\oplus 111^d$ for some $a,b,c,d\in\mathbb{Z}_{\geq 0}$.  By direct calculation, $\JFg(V)=((c+d),(\max(b,c)+a+d, \min(b,c)), (b+d))$. We obtain the corresponding reverse plane partition shown in Figure~\ref{fig:introRPP}.

\begin{figure}[!htbp]
\[\raisebox{.65in}{$010^a\oplus 011^b\oplus 110^c\oplus 111^d \hspace{.2in}\longleftrightarrow$\hspace{.2in}} 
\begin{tikzpicture}[scale=1.5]
\node (1) at (1,0) {{$b+d$}};
\node (3) at (2,1) {{$\min(b,c)$}};
\node (2) at (0,1) {{$\max(b,c)+a+d$}};
\node (4) at (1,2) {{$c+d$}};
\draw[->] (1)--(3);
\draw[->] (2)--(1);
\draw[->] (2)--(4);
\draw[->] (4)--(3);
\end{tikzpicture}\]
\caption{The correspondence between isomorphism classes of representations of $Q=1\rightarrow 2 \leftarrow 3$ belonging to $\mathcal{C}_{Q,2}$ and reverse plane partitions of $\textsf{P}_{Q,2}$.}
\label{fig:introRPP}
\end{figure}
\end{example}
\subsection{Periodicity}
For this section, we assume that the vertices of $Q$ are labeled with the numbers $1,\ldots, n$ so that the arrows of $Q$ are oriented from vertices with larger labels to vertices with smaller labels.

Note that $t_{\sf x}$ and $t_{\sf y}$ commute unless $\sf x$ and $\sf y$ are related by a cover.  
For $\pos$, and $i\in Q_0$, the elements of $\tau$-orbit corresponding to $i$ are never related by a cover, so we can define ${T}_i=\prod_{\sf x} t_{\sf x}$, where we toggle at each element {$\sf x$ of the $\tau$-orbit containing the projective at $i$}, without worrying about the order in which the composition is taken.

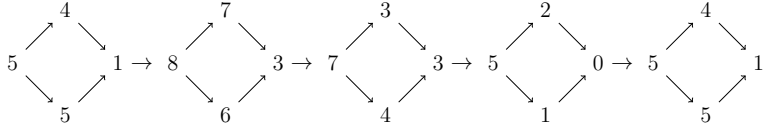
\begin{figure}
    \centering
    \scalebox{.7}{
    \begin{tikzpicture}
    \node (left) at (-1,1) {5};
    \node (top) at (0,2) {4};
    \node (bot) at (0,0) {5};
    \node (right) at (1,1) {1};
    \draw[->] (left)--(bot);
    \draw[->] (left)--(top);
    \draw[->] (top)--(right);
    \draw[->] (bot)--(right);
    \end{tikzpicture}\raisebox{.45in}{$\rightarrow$}
\begin{tikzpicture}
    \node (left) at (-1,1) {8};
    \node (top) at (0,2) {7};
    \node (bot) at (0,0) {6};
    \node (right) at (1,1) {3};
    \draw[->] (left)--(bot);
    \draw[->] (left)--(top);
    \draw[->] (top)--(right);
    \draw[->] (bot)--(right);
    \end{tikzpicture}\raisebox{.45in}{$\rightarrow$}
    \begin{tikzpicture}
    \node (left) at (-1,1) {7};
    \node (top) at (0,2) {3};
    \node (bot) at (0,0) {4};
    \node (right) at (1,1) {3};
    \draw[->] (left)--(bot);
    \draw[->] (left)--(top);
    \draw[->] (top)--(right);
    \draw[->] (bot)--(right);
    \end{tikzpicture}\raisebox{.45in}{$\rightarrow$}
    \begin{tikzpicture}
    \node (left) at (-1,1) {5};
    \node (top) at (0,2) {2};
    \node (bot) at (0,0) {1};
    \node (right) at (1,1) {0};
    \draw[->] (left)--(bot);
    \draw[->] (left)--(top);
    \draw[->] (top)--(right);
    \draw[->] (bot)--(right);
    \end{tikzpicture}\raisebox{.45in}{$\rightarrow$}
    \begin{tikzpicture}
    \node (left) at (-1,1) {5};
    \node (top) at (0,2) {4};
    \node (bot) at (0,0) {5};
    \node (right) at (1,1) {1};
    \draw[->] (left)--(bot);
    \draw[->] (left)--(top);
    \draw[->] (top)--(right);
    \draw[->] (bot)--(right);
    \end{tikzpicture}}
    \caption{An example of periodicity, where $N=8$. In each step, we first toggle the elements in the middle $\tau$-orbit, and then we toggle the top and bottom $\tau$-orbits.}
    \label{fig:periodicity}
\end{figure}

Define $\pro_Q={T}_n\circ \dots \circ {T}_1$.  Define $h$ to be the Coxeter number of $Q$: by definition, this is the order of the product of the simple reflections in the Coxeter group, or, equivalently, the largest degree of the root system. Figure~\ref{fig:periodicity} illustrates the following theorem.

\begin{theorem}\label{intro-period-N} For $Q$ a Dynkin quiver and $m$ a minuscule vertex, $\pro_Q^h$ is the identity transformation on RPPs on
$\pos$ with entries in $[0,N]$.
\end{theorem}





\section{Type \emph{A} examples}
Suppose we start with the type $A$ quiver $Q$ shown in Figure~\ref{fig:HGexample}, where 3 is the chosen minuscule vertex. The corresponding AR quiver is also shown in the figure. The type $A$ minuscule poset $\pos$ associated with vertex 3 is shown in black in the figure. 

We denote a representation $M$ as a labeling of the poset $\pos$, where the label at a vertex denotes how many copies of the corresponding indecomposable are in $M$. The representation $M$ in Figure~\ref{fig:HGexample} contains 4 copies of the indecomposable with dimension vector 11100 and 3 copies of the indecomposable with dimension vector 01100. 

\begin{figure}[!htbp]
\centering
\scalebox{.8}{
\begin{tikzpicture}[scale=1]
\node (11100) at (1,0) {$11100_{13}$};
\node (00010) at (3,0) {\textcolor{gray}{$00010_8$}};
\node (00001) at (5,0) {\textcolor{gray}{$00001_3$}}; 
\node (01100) at (0,1) {$01100_{14}$};
\node (11110) at (2,1) {$11110_9$};
\node (00011) at (4,1) {\textcolor{gray}{$00011_4$}};
\node (01110) at (1,2) {$01110_{10}$};
\node (11111) at (3,2) {$11111_5$};
\node (01111) at (2,3) {$01111_7$};
\node (11000) at (4,3) {\textcolor{gray}{$11000_2$}};
\node (01000) at (3,4) {\textcolor{gray}{$01000_6$}};
\node (00100) at (-1,2) {$00100_{15}$};
\node (00110) at (0,3) {$00110_{12}$};
\node (00111) at (1,4) {$00111_{11}$};
\node (10000) at (5,4) {\textcolor{gray}{$10000_1$}};
\draw[->] (00100)--(01100);
\draw[->] (01100)--(11100);
\draw[->] (00100)--(00110);
\draw[->] (00110)--(00111);
\draw[->] (00110)--(01110);
\draw[->] (01100)--(01110);
\draw[->] (11100)--(11110);
\draw[->] (11110)--(11111);
\draw[->] (11111)--(11000);
\draw[->] (11000)--(10000);
\draw[->] (01110)--(11110);
\draw[->] (11110)--(00010);
\draw[->] (00111)--(01111);
\draw[->] (01111)--(11111);
\draw[->] (11111)--(00011);
\draw[->] (00011)--(00001);
\draw[->] (00010)--(00011);
\draw[->] (01110)--(01111);
\draw[->] (01111)--(01000);
\draw[->] (01000)--(11000);
\end{tikzpicture}
\raisebox{1in}{$M=$
\raisebox{-.65in}{\begin{tikzpicture}[scale=.7]
\node (11100) at (1,0) {4};
\node (01100) at (0,1) {3};
\node (11110) at (2,1) {0};
\node (01110) at (1,2) {1};
\node (11111) at (3,2) {1};
\node (01111) at (2,3) {0};
\node (00100) at (-1,2) {0};
\node (00110) at (0,3) {1};
\node (00111) at (1,4) {2};
\draw[->] (00100)--(01100);
\draw[->] (01100)--(11100);
\draw[->] (00100)--(00110);
\draw[->] (00110)--(00111);
\draw[->] (00110)--(01110);
\draw[->] (01100)--(01110);
\draw[->] (11100)--(11110);
\draw[->] (11110)--(11111);
\draw[->] (01110)--(11110);
\draw[->] (00111)--(01111);
\draw[->] (01111)--(11111);
\draw[->] (01110)--(01111);
\end{tikzpicture}}}}
\scalebox{.8}{
\begin{ytableau}
1 &\textcolor{lightgray}{0}&\textcolor{lightgray}{0} \\
\textcolor{lightgray}{0} & \textcolor{lightgray}{0}& \textcolor{lightgray}{0}\\
\textcolor{lightgray}{0} & \textcolor{lightgray}{0}& \textcolor{lightgray}{0}
\end{ytableau}\hspace{.1in}
\begin{ytableau}
1 & 1 & \textcolor{lightgray}{0}\\
\textcolor{lightgray}{0} & \textcolor{lightgray}{0}& \textcolor{lightgray}{0}\\
\textcolor{lightgray}{0} & \textcolor{lightgray}{0}& \textcolor{lightgray}{0}
\end{ytableau}\hspace{.1in}
\begin{ytableau}
1 & 1 & \textcolor{lightgray}{0}\\
1 & \textcolor{lightgray}{0}& \textcolor{lightgray}{0}\\
\textcolor{lightgray}{0} & \textcolor{lightgray}{0}& \textcolor{lightgray}{0}
\end{ytableau}\hspace{.1in}
\begin{ytableau}
0 & 1 & \textcolor{lightgray}{0}\\
1 & 2 & \textcolor{lightgray}{0}\\
\textcolor{lightgray}{0} & \textcolor{lightgray}{0}& \textcolor{lightgray}{0}
\end{ytableau}\hspace{.1in}
\begin{ytableau}
0 & 1 & 3\\
1 & 2 & 4\\
\textcolor{lightgray}{0} & \textcolor{lightgray}{0}& \textcolor{lightgray}{0}
\end{ytableau}\hspace{.1in}
\begin{ytableau}
0 & 1 & 3\\
1 & 2 & 4\\
5 &\textcolor{lightgray}{0} &\textcolor{lightgray}{0} 
\end{ytableau}\hspace{.1in}
\begin{ytableau}
0 & 1 & 3\\
1 & 2 & 4\\
5 & 8 & \textcolor{lightgray}{0}
\end{ytableau}\hspace{.1in}
\begin{ytableau}
1 & 1 & 3\\
1 & 3 & 4\\
5 & 8 & 8
\end{ytableau}
\raisebox{-.15in}{ $=\RSK$}
}
\caption{The top left image is the AR quiver associated with the quiver $Q=1\rightarrow 2\rightarrow \mathbf{3}\leftarrow 4 \leftarrow 5$ with chosen minuscule vertex 3. The dimension vectors with support on vertex 3 are black, while the others are gray. The top right image represents a representation $M$ whose indecomposable summands all have support on vertex 3. The images below show the steps in computing $\RSK$, shown on the corresponding Young diagram.}
\label{fig:Pakexample}
\end{figure}

The order we use to compute $\RSK$ is indicated in the subscripts on the AR quiver. The procedure is shown step by step in the 12 fillings in Figure~\ref{fig:HGexample}. In future sections, it will be useful to realize the resulting reverse plane partition as a reverse plane partition for the Young diagram of shape $(3,3,3)$, as shown.

Figure~\ref{fig:Pakexample} shows another example using a different orientation $Q$ and a representation that assigns the same multiplicities to the representations with the same dimension vectors as above. In this figure, we show how to carry out the algorithm by identifying the intermediate fillings with fillings of a Young diagram. 

A reader who is familiar with the known bijections between multisets of rim hooks of a Young diagram and reverse plane partitions of the same Young diagram can check that the example in Figure~\ref{fig:HGexample} agrees with the well-known Hillman--Grassl correspondence and the example in Figure~\ref{fig:Pakexample} agrees with the bijection first described by Pak \cite{pak2001hook} and later by Sulzgruber \cite{sulzgruberfull}  and Hopkins \cite{hopkins2014rsk}. In fact, this is not a coincidence, and we make the correspondence precise in Sections~\ref{sec:HG} and ~\ref{sec:Pak}.

\section{The Hillman--Grassl Correspondence}\label{sec:HG}

A \textit{rim hook} of the Young diagram $\lambda$ is a connected strip of border boxes in $\lambda$ such that the result of removing these boxes is again a Young diagram. The Hillman--Grassl correspondence \cite{hillman1976reverse} is a bijection between multisets of rim hooks of $\lambda$ and the set of reverse plane partitions of $\lambda$. See, for example, \cite{stanley1999enumerative}. 

\begin{figure}[!htbp]
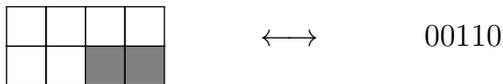

\centering
\begin{ytableau}  $ $ & & & \\
$ $ & & *(gray) & *(gray)
\end{ytableau}\hspace{.5in}$\longleftrightarrow$\hspace{.5in} 00110
\caption{A rim hook and the corresponding dimension vector.}
\label{fig:hooktodimvector}
\end{figure}

Let $M$ be a representation in $\cat$, where $Q$ is a type $A$ Dynkin quiver. Then as shown in Figure~\ref{fig:HGexample},  $\RSK$ can be viewed as a reverse plane partition of the Young diagram of shape $m^{(n+1-m)}$. We may identify the indecomposable summands of $M$ as rim hooks of $m^{(n+1-m)}$ by reading through the southeast border of $m^{(n+1-m)}$ from southwest to northeast and including a box in the rim hook exactly when the corresponding entry in the dimension vector of the indecomposable is 1. See Figure~\ref{fig:hooktodimvector}. In this way, we may identify $M$ with a multiset of rim hooks of $m^{(n+1-m)}$. Given a multiset $M$ of rim hooks of $m^{(n+1-m)}$ (i.e., an $M$ in $\cat$), let $HG(M)$ denote the reverse plane partition of shape $m^{(n+1-m)}$ obtained using the Hillman--Grassl correspondence. 

\begin{theorem}\label{thm:HG} Let $Q$ be a type $A_n$ quiver with chosen minuscule vertex $m$ and the following orientation.
\[Q=1\leftarrow 2 \leftarrow \cdots \leftarrow n-1 \leftarrow n\]
Then for any $M\in \cat$, $\RSK=HG(M)$.
\end{theorem}

\begin{remark}\label{remark:HGrectangle}
We can, in fact, recover the Hillman--Grassl correspondence for any shape $\lambda$---not just for rectangles---using $\rho_{Q,m}(-)$. 

\end{remark}

\section{The Pak correspondence}\label{sec:Pak}

In \cite{pak2001hook}, Pak gives a different bijection between multisets of rim hooks of a partition $\lambda$ and reverse plane partitions of $\lambda$, which involves toggles. We again use the identification between rim hooks of a rectangle and dimension vectors of indecomposable representations in $\cat$ from Section~\ref{sec:HG}. Let $Pak(M)$ denote the reverse plane partition obtained from Pak's correspondence using the multiset of rim hooks determined by $M$. The following result is illustrated in Figure~\ref{fig:Pakexample}.

\begin{theorem}\label{thm:Pak} Let $Q$ be a type $A_n$ quiver with chosen minuscule vertex $m$ and the following orientation.
\[Q=1\rightarrow \cdots \rightarrow m \leftarrow \cdots \leftarrow n\]
Then for any $M\in \cat$, $\RSK=Pak(M)$.
\end{theorem}

\begin{remark}
We can, in fact, recover the Pak correspondence for any shape $\lambda$ using $\rho_{Q,m}(-)$.
\end{remark}

Fix $n$ and $m\in\{1,\ldots,n\}$. For each type $A_n$ Dynkin quiver, we obtain a bijection $\rho_{Q,m}(-)$ that we can think of as a map between multisets of rim hooks of $m^{(n+1-m)}$ and reverse plane partitions of $m^{(n+1-m)}$. We have shown that for two particular orientations, we recover the Hillman-Grassl and Pak correspondences. For other orientations, we obtain different bijections. 
We demonstrate this with an example. 

We begin with the quiver 
\[Q=1\leftarrow 2\rightarrow \mathbf{3}\leftarrow 4\rightarrow 5\] with minuscule vertex 3 and use the following multiset of rim hooks of the 3-by-3 rectangle.
\[01100,01110,11111,11111,00110,00110\]
This multiset of rim hooks determines the representation $M$ of the quiver $Q$ shown below. We then obtain $\RSK$ as shown. We leave it to the reader to check that both Hillman--Grassl and the Pak correspondences applied to this same multiset of rim hooks yield a reverse plane partition where the maximal element is labeled with 0.
\begin{center}
\scalebox{.8}{
\begin{tikzpicture}[scale=1]
\node (00110) at (1,0) {$00110_{9}$};
\node (11000) at (3,0) {\textcolor{gray}{$11000_4$}};
\node (00111) at (0,1) {$00111_{11}$};
\node (11110) at (2,1) {$11110_6$};
\node (01000) at (4,1) {\textcolor{gray}{$01000_1$}};
\node (11111) at (1,2) {$11111_{8}$};
\node (01110) at (3,2) {$01110_3$};
\node (01111) at (2,3) {$01111_7$};
\node (00010) at (4,3) {\textcolor{gray}{$00010_2$}};
\node (00011) at (3,4) {\textcolor{gray}{$00011_5$}};
\node (00100) at (-1,2) {$00100_{13}$};
\node (11100) at (0,3) {$11100_{12}$};
\node (01100) at (1,4) {$01100_{10}$};
\node (00001) at (-1,0) {\textcolor{gray}{$00001_{14}$}};
\node (10000) at (-1,4) {\textcolor{gray}{$10000_{15}$}};
\draw[->] (00001)--(00111);
\draw[->] (00100)--(00111);
\draw[->] (00100)--(11100);
\draw[->] (10000)--(11100);
\draw[->] (00111)--(00110);
\draw[->] (00111)--(11111);
\draw[->] (11100)--(11111);
\draw[->] (11100)--(01100);
\draw[->] (00110)--(11110);
\draw[->] (11111)--(11110);
\draw[->] (11111)--(01111);
\draw[->] (01100)--(01111);
\draw[->] (11110)--(11000);
\draw[->] (11110)--(01110);
\draw[->] (01111)--(01110);
\draw[->] (01111)--(00011);
\draw[->] (00011)--(00010);
\draw[->] (01110)--(00010);
\draw[->] (01110)--(01000);
\draw[->] (11000)--(01000);
\end{tikzpicture}
\raisebox{.85in}{$M=$}
\raisebox{.23in}{
\begin{tikzpicture}[scale=.7]
\node (00110) at (1,0) {$2$};
\node (00111) at (0,1) {$0$};
\node (11110) at (2,1) {$0$};
\node (11111) at (1,2) {$2$};
\node (01110) at (3,2) {$1$};
\node (01111) at (2,3) {$0$};
\node (00100) at (-1,2) {$0$};
\node (11100) at (0,3) {$0$};
\node (01100) at (1,4) {$1$};
\draw[->] (00100)--(00111);
\draw[->] (00100)--(11100);
\draw[->] (00111)--(00110);
\draw[->] (00111)--(11111);
\draw[->] (11100)--(11111);
\draw[->] (11100)--(01100);
\draw[->] (00110)--(11110);
\draw[->] (11111)--(11110);
\draw[->] (11111)--(01111);
\draw[->] (01100)--(01111);
\draw[->] (11110)--(01110);
\draw[->] (01111)--(01110);
\end{tikzpicture}}
\raisebox{.85in}{$\RSK=$}
\raisebox{.23in}{
\begin{tikzpicture}[scale=.7]
\node (00110) at (1,0) {$2$};
\node (00111) at (0,1) {$3$};
\node (11110) at (2,1) {$2$};
\node (11111) at (1,2) {$2$};
\node (01110) at (3,2) {$1$};
\node (01111) at (2,3) {$1$};
\node (00100) at (-1,2) {$3$};
\node (11100) at (0,3) {$3$};
\node (01100) at (1,4) {$2$};
\draw[->] (00100)--(00111);
\draw[->] (00100)--(11100);
\draw[->] (00111)--(00110);
\draw[->] (00111)--(11111);
\draw[->] (11100)--(11111);
\draw[->] (11100)--(01100);
\draw[->] (00110)--(11110);
\draw[->] (11111)--(11110);
\draw[->] (11111)--(01111);
\draw[->] (01100)--(01111);
\draw[->] (11110)--(01110);
\draw[->] (01111)--(01110);
\end{tikzpicture}}}
\end{center}

\bibliography{sample}{}
\bibliographystyle{plain}


\end{document}